\documentclass[10pt]{amsart}
\usepackage{amsmath,amssymb}
\newcommand{\doublespace}
   {\addtolength{\baselineskip}{0.15\baselineskip}}

\newtheorem{pdef}{Definition}[section] %
\newtheorem{thm}[pdef]{Theorem}        
\newtheorem{cor}[pdef]{Corollary}
\newtheorem{lem}[pdef]{Lemma}
\newtheorem{prop}[pdef]{Proposition}

\newcounter{equationnumber}
\renewcommand{\theequation}{\thesection.\arabic{equation}}
\def\mathletters{
    \addtocounter{equation}{1}
    \edef\@currentlabel{\theequation}
    \setcounter{equationnumber}{\value{equation}}
    \setcounter{equation}{0}
    \edef\theequation{\@currentlabel\noexpand\alph{equation}}
    }

\title{Supports of Measures in a free additive convolution semigroup}
\author{Hao-Wei Huang}
\address{Department of Mathematics, Indiana University, 831 East 3rd Street,
Bloomington, IN 47405} \email{huang39@indiana.edu}

\begin{document}
\maketitle \doublespace \pagestyle{myheadings} \thispagestyle{plain}
\markboth{   }{ }

\begin{abstract} In this paper, we study
the supports of measures in the free additive convolution semigroup
$\{\mu^{\boxplus t}:t>1\}$, where $\mu$ is a Borel probability
measure on $\mathbb{R}$. We give a formula for the density of the
absolutely continuous part of $\mu^{\boxplus t}$ and use this
formula to obtain certain regularizing properties of $\mu^{\boxplus
t}$. We show that the number $n(t)$ of the components in the support
of $\mu^{\boxplus t}$ is a decreasing function of $t$ and give
equivalent conditions so that $n(t)=1$ for sufficiently large $t$.
Moreover, a measure $\mu$ so that $\mu^{\boxplus t}$ has infinitely
many components in the support for all $t>1$ is given.
\end{abstract}

\section{Introduction}

Given two Borel probability measures $\mu$ and $\nu$ on
$\mathbb{R}$, denote by $\mu\boxplus\nu$ their free additive
convolution, which is the probability distribution of $X+Y$, where
$X$ and $Y$ are free selfadjoint random variables with distributions
$\mu$ and $\nu$, respectively. Free convolution is a binary
operation on the set of Borel probability measures on $\mathbb{R}$,
which is analogue of classical convolutions. We refer to [19] for a
systematic exposition of this subject and [5] for analytic methods
for the computations of free convolutions.

One of the important properties for free convolution is
subordination. That is, the Cauchy transform (see section 2 for
definition) $G_{\mu\boxplus\nu}$ is subordinated to $G_\mu$, in the
sense that $G_{\mu\boxplus\nu}=G_\mu\circ\omega$ for some analytic
self-mapping $\omega$ of the complex upper half-plane. We refer the
reader to [7,9, and 17] for details. One of the applications of
subordination functions is to study the regularity for free
convolution. Free convolution is a highly nonlinear operation,
however, it has a stronger regularizing effect than classical
convolution. Systematic study of regularity for free convolution by
means of subordination functions can be found in [6,7,8, and 17]. In
particular, the atoms of free convolutions are determined. Moreover,
if $\gamma$ is the semicircular distribution then the density of
$\mu\boxplus\gamma$ is continuous on $\mathbb{R}$ and analytic
wherever it is positive.

For $n\in\mathbb{N}$, the $n$-fold free convolution
$\mu\boxplus\cdots\boxplus\mu$ is denoted by $\mu^{\boxplus n}$. It
is known that the discrete semigroup $\{\mu^{\boxplus
n}:n\in\mathbb{N}\}$ can be embedded in a continuous family
$\{\mu^{\boxplus t}:t\geq1\}$ such that $\mu^{\boxplus
t_1}\boxplus\mu^{\boxplus t_2}=\mu^{\boxplus(t_1+t_2)}$ for all
$t_1,t_2\geq1$. This was first proved by Bercovici and Voiculescu in
[6] for measures $\mu$ with compact support and for large values of
$t$. Later, in [14] Nica and Speicher generalized this result to
compactly supported measures and $t>1$ by exhibiting explicit random
variables with the corresponding distributions whose calculations of
the distribution are based on moments and combinatorics. The
existence of subordination functions $\omega_t$ in this context and
related consequences for the regularity for the measures
$\mu^{\boxplus t}$ were thoroughly studied by Belinschi and
Bercovici in [1] and [2]. More precisely, for $t>1$ let
$\mu^{\boxplus t}$ be the measure satisfying the requirement
$F_{\mu^{\boxplus t}}=F_\mu\circ\omega_t$, where $F_{\mu^{\boxplus
t}}=1/G_{\mu^{\boxplus t}}$ and $F_\mu=1/G_\mu$ are reciprocal
Cauchy transforms. This gives an alternate proof of the existence of
$\mu^{\boxplus t}$. With the help of $\omega_t$, it was shown that
the free convolution power $\mu^{\boxplus t}$ of $\mu$ has no
singular continuous part and its density in the absolutely
continuous part is locally analytic.

In this present paper, we use different points of views to obtain
some properties of $F_\mu$ and analyze the relation between the
supports $\mathrm{supp}(\mu)$ and $\mathrm{supp}(\mu^{\boxplus t})$
of $\mu$ and $\mu^{\boxplus t}$, $t>1$, respectively. By means of
the methods developed in this paper, we give a complete description
about how $\mathrm{supp}(\mu^{\boxplus t})$ changes when $t$
increases. Motivated by the free central limit theorem (see
[6,10,11,12, 15, and 20] for details), one of the purposes of this
paper is to study the number $n(t)$ of components in the support of
$\mu^{\boxplus t}$. We show that $n(t)$ is a decreasing function of
$t$ for $t>1$ and that $n(t)=1$ for sufficiently large $t$ if and
only if $n(t')$ is finite for some $t'>1$. We also construct a
measure $\mu$ so that $n(t)$ is infinite for all $t>1$.

The paper is organized as follows. Section 2 contains some
definitions and basic propositions from free probability theory.
Section 3 provides some properties of $F_\mu$ which are related to
the subordination function $\omega_t$, and an explicit formula for
the density of the absolutely continuous part of $\mu^{\boxplus t}$.
Section 4 contains complete investigations on the supports of the
measures $\mu^{\boxplus t}$.

\section{Preliminary}

For any complex $z$ in $\mathbb{C}$, let $\Re z$ and $\Im z$ be the
real and imaginary parts of $z$, respectively, and denote by
$\mathbb{C}^+=\{z\in\mathbb{C}:\Im z>0\}$ the complex upper
half-plane. For any Borel probability measure $\mu$ on $\mathbb{R}$,
define its Cauchy transform as
\[G_\mu(z)=\int_\mathbb{R}\frac{d\mu(s)}{z-s},\;\;\;z\in\mathbb{C}^+.\] As shown in [5], the composition
inverse $G_\mu^{-1}$ of $G_\mu$ is defined in an appropriate Stolz
angle
\[D_{\eta,\epsilon}=\{x+iy:|x|<-\eta y,\;0<-y<\epsilon\}\]
of zero in the lower half-plane. The Voiculescu transform
$\mathcal{R}_\mu$ of the measure $\mu$ defined as
\[\mathcal{R}_\mu(z)=G_\mu^{-1}(z)-z^{-1},\;\;\;z\in D_{\eta,\epsilon},\]
is the linearizing transform for the operation $\boxplus$. That is,
the identity
\[\mathcal{R}_{\mu\boxplus\nu}=\mathcal{R}_\mu+\mathcal{R}_\nu\]
holds on a domain where these three functions are defined. The
measure $\mu$ can be recovered from the Cauchy transform $G_\mu$ as
the weak$^*$-limit of the measures
\[d\mu_\epsilon(s)=-\frac{1}{\pi}\Im G_\mu(s+i\epsilon)\;ds\] as $\epsilon\to0^+$, which
we refer to as the Stieltjes inversion formula. Particularly, if the
function $\Im G_\mu(x)$ is continuous at the point $x\in\mathbb{R}$
then in the Lebesgue decomposition the density of the absolutely
continuous part of $\mu$ at $x$ is given by $-\Im G_\mu(x)/\pi$.
Therefore, this inversion formula gives a way to extract the density
function from the Cauchy transform.

The reciprocal Cauchy transform $F_\mu=1/G_\mu$ is an analytic
self-mapping of $\mathbb{C}^+$. By Nevanlinna representation, it can
be expressed as
\[F_\mu(z)=a+z+\int_\mathbb{R}\frac{1+sz}{s-z}\;d\rho(s),\leqno{(2.1)}\] where $a=\Re
F_\mu(i)$ and $\rho$ is some finite positive Borel measure on
$\mathbb{R}$ uniquely determined by $F_\mu$. Indeed, $\rho$ is the
weak$^*$-limit as $\epsilon\to0^+$ of the measures
\[d\rho_\epsilon(s)=\frac{\Im F_\mu(s+i\epsilon)}{\pi(s^2+1)}\;ds.\leqno{(2.2)}\] It was
shown in [13] that $F_\mu$ can be expressed in terms of a Cauchy
transform of some measure if $\mu$ has finite mean and unit
variance. More precisely, a measure $\mu$ satisfies
\[\int_\mathbb{R}s\;d\mu(s)=0\;\;\;\;\mathrm{and}\;\;\;\;\int_\mathbb{R}s^2\;d\mu(s)=1\] if and
only if \[F_\mu(z)=z-G_\nu(z),\;\;z\in\mathbb{C}^+,\] where $\nu$ is
some Borel probability measure on $\mathbb{R}$.

Denote by $\mathop{z\to\alpha}\limits_{\sphericalangle}$ if
$z\to\alpha$ nontangentially to $\mathbb{R}$, i.e., $(\Re
z-\alpha)/\Im z$ stays bounded. A useful criterion for locating an
atom $\alpha$ of $\mu$ is the limit
\[\lim_{\mathop{z\to\alpha}\limits_{\sphericalangle}}
(z-\alpha)G_\mu(z)=\mu(\{\alpha\}).\leqno{(2.3)}\] An equivalent
statement of the above limit is that a point $\alpha\in\mathbb{R}$
is an atom of $\mu$ if and only if $F_\mu(\alpha)=0$ and the
Julia-Carath\'{e}odoty derivative $F_\mu'(\alpha)$ (which is the
limit of
\[\frac{F_\mu(z)-F_\mu(\alpha)}{z-\alpha}\] as $z\to\alpha$ nontangentially and $z\in\mathbb{C}^+$)
is finite, in which case
\[F_\mu'(\alpha)=\frac{1}{\mu(\{\alpha\})}.\]
Similarly, there is an effective way to identify the atoms of the
measure $\rho$ in the Nevanlinna representation (2.1) for $F_\mu$.
The following lemma is well known, however, we have no recollection
of seeing it stated explicitly. We provide a full proof without
claiming any paternity.

\begin{lem} If $\alpha\in\mathbb{R}$ then
\[\lim_{\mathop{z\to\alpha}\limits_{\sphericalangle}}(z-\alpha)F_\mu(z)=-(\alpha^2+1)\rho(\{\alpha\}).\]
\end{lem}

\begin{pf} By considering the positive measure $\rho-\rho(\{\alpha\})\delta_0$, we
may assume that $\alpha$ is not an atom of $\rho$, in which case we
must show that $(z-\alpha)F_\mu(z)\to0$ as $z\to\alpha$
nontangentially. Write $z=x+iy$ and suppose
\[\frac{|x-\alpha|}{y}\leq
M\;\;\;\;\mathrm{and}\;\;\;\;|z-\alpha|<1,\] where $0<M<\infty$. Let
$N=2(|\alpha|+1)$. Note that if $s\geq N$ then
\[\left|\frac{s}{s-z}\right|=\left|\frac{\frac{s}{y}}{\frac{s-x}{y}-i}\right|\leq
\frac{s}{s-x}\leq\frac{s}{s-(\alpha+1)}\leq2,\] whence
\[|z-\alpha|\left|\frac{sz}{s-z}\right|\leq2(|\alpha|+1).\]
Similarly, the above inequality also holds for $s\leq-N$. For
$s\in(-N,N)$, we have
\[|z-\alpha|\left|\frac{sz}{s-z}\right|\leq
\frac{|z-\alpha|}{y}|sz|\leq N(|\alpha|+1)\sqrt{M^2+1}.\] Moreover,
for any $s\in\mathbb{R}$ we have
\[\frac{|z-\alpha|}{|s-z|}=\frac{\left|\frac{x-\alpha}{y}+i\right|}{\left|\frac{x-s}{y}+i\right|}\leq
\left|\frac{x-\alpha}{y}+i\right|\leq\sqrt{M^2+1}.\] Since
\[\left|(z-\alpha)\int_\mathbb{R}\frac{1+sz}{s-z}\;d\rho(s)\right|\leq
\int_\mathbb{R}\frac{|z-\alpha|+|sz||z-\alpha|}{|s-z|}\;d\rho(s),\]
and as $z\to\alpha$ nontangentially the integrand converges to zero
$\rho$-a.e. and stays bounded by the above discussions, it is easy
to see the desired result holds.
\end{pf} \qed

As in the introduction, the collection of the $t$-th convolution
power $\mu^{\boxplus t}$ of $\mu$ forms a semigroup $\{\mu^{\boxplus
t}:t>1\}$ which interpolates the discrete semigroup $\{\mu^{\boxplus
n}:n\in\mathbb{N}\}$. The construction of $\mu^{\boxplus t}$ for
arbitrary $\mu$ by using subordination functions and analytic
methods is introduced in following two paragraphs.

For any $t>1$, consider the function
\[H_t(z)=tz-(t-1)F_\mu(z),\;\;\;z\in\mathbb{C}^+.\leqno{(2.4)}\] Then we have $\Im H_t(z)\leq\Im
z$ for $z\in\mathbb{C}^+$ and
\[\lim_{y\to\infty}\frac{H_t(iy)}{iy}=1.\] Moreover, there exists
a continuous function
$\omega_t:\overline{\mathbb{C}^+}\to\overline{\mathbb{C}^+}$ such
that $\omega_t(\mathbb{C}^+)\subset\mathbb{C}^+$,
$\omega_t|\mathbb{C}^+$ is analytic, and $H_t(\omega_t(z))=z$ for
$z\in\mathbb{C}^+$. If the set $\{z\in\mathbb{C}^+:\Im H_t(z)>0\}$
is denoted by $\Omega_t$ then
\[\Omega_t=\omega_t(\mathbb{C}^+)\]
is a simply connected set whose boundary $\omega_t(\mathbb{R})$ is a
simple curve. Then the equation
$\omega_t(H_t(\omega_t(z)))=\omega_t(z)$, $z\in\mathbb{C}^+$, shows
that the equation $\omega_t(H_t(z))=z$ holds for $z\in\Omega_t$, and
hence the inverse of $\omega_t$ is $H_t|\Omega_t$. Moreover, if
$x\in\mathbb{R}$ and $\Im\omega_t(x)>0$ then $\omega_t$ can be
continued analytically to a neighborhood of $x$.

Let $\mu^{\boxplus t}$ be the measure defined by the requirement
\[F_{\mu^{\boxplus
t}}(z)=F_\mu(\omega_t(z)),\;\;\;\;z\in\mathbb{C}^+.\leqno{(2.5)}\]
It turns out that the measure $\mu^{\boxplus t}$ can be
characterized in terms of Voiculescu transform, i.e., the measure
$\mu^{\boxplus t}$ obtained in this way is the unique measure
satisfying
\[\mathcal{R}_{\mu^{\boxplus t}}(z)=t\mathcal{R}_\mu(z),\] where $z$
is in the common domain of these two functions. By using the
equation $H_t(\omega_t(z))=z$ we can rewrite (2.5) as
\[F_{\mu^{\boxplus
t}}(z)=\frac{t\omega_t(z)-z}{t-1},\;\;\;z\in\mathbb{C}^+,\] which
shows, in particular, that $F_{\mu^{\boxplus t}}$ extends
continuously to $\overline{\mathbb{C}^+}$. Also, by the equation
$\omega_t(H_t(z))=z$, $z\in\Omega_t$, we have
\[F_{\mu^{\boxplus t}}(H_t(z))=F_\mu(z),\;\;\;z\in\Omega_t.\]

Regularity properties of $\mu^{\boxplus t}$ have been studied
thoroughly in [1] and [2]. We now review some of these results which
are relevant to our investigation on the number of the components in
the support of $\mu^{\boxplus t}$. For any $t>1$, the singular part
of the measure $\mu^{\boxplus t}$ with respect to Lebesgue measure
is purely atomic and $\mu^{\boxplus t}$ has an atom $\alpha$ if and
only if $\mu(\{t\alpha\})>1-t^{-1}$ in which case $\mu^{\boxplus
t}(\{\alpha\})=t\mu(\{\alpha\})-(t-1)$. A weaker statement, i.e., a
point $\alpha$ satisfying $F_{\mu^{\boxplus t}}(t\alpha)=0$ if and
only if $\mu(\{\alpha\})\geq1-t^{-1}$, also holds. Hence for any
measure $\mu$ which is not a point mass, the atoms of $\mu^{\boxplus
t}$ disappear for sufficiently large $t$, and hence $\mu^{\boxplus
t}$ is absolutely continuous for such $t$.

\section{Decreasing number of components in the support of $\mu^{\boxplus t}$}

Throughout this section, $t$ is a parameter with $t>1$ and for any
given and fixed Borel probability measure $\mu$ on $\mathbb{R}$,
$\rho$ is the unique finite positive Borel measure in the Nevanlinna
representation (2.1) of $F_\mu$. Our analysis of the support of
$\mu^{\boxplus t}$ will be based on the following function
$g:\mathbb{R}\to\mathbb{R}_+\cup\{\infty\}$ defined as
\[g(x)=\int_\mathbb{R}\frac{s^2+1}{(s-x)^2}\;d\rho(s),\;\;\;x\in\mathbb{R}.\]
Associated with the function $g$ and the parameter $t$ are the
following sets
\[V_t^+=\left\{x\in\mathbb{R}:g(x)>\frac{1}{t-1}\right\},\]
\[V_t=\left\{x\in\mathbb{R}:g(x)=\frac{1}{t-1}\right\},\]
\[V_t^-=\left\{x\in\mathbb{R}:g(x)<\frac{1}{t-1}\right\},\]
and the function
\[f_t(x)=\inf\left\{y\geq0:\int_\mathbb{R}\frac{s^2+1}{(x-s)^2+y^2}\;d\rho(s)\leq\frac{1}{t-1}\right\}\]
which will play important roles in the analysis.

\begin{lem} If $V_t\cup V_t^-$ contains an open interval $I$
then $\rho(I)=0$ and $g$ is strictly convex on $I$.
\end{lem}

\begin{pf} Let $(a,b)$ be any interval contained in $I$. Then for
any $x\in(a,b)$ we have
\[\frac{1}{t-1}\geq\int_a^b\frac{s^2+1}{(s-x)^2}\;d\rho(s)\geq
\int_a^b\frac{1}{(b-a)^2}\;d\rho(s)=\frac{\rho((a,b))}{(b-a)^2},\]
from which we obtain that $\rho((a,b))=0$. Since the second order
derivative of $g$ is positive on $(a,b)$, the second assertion
follows. This completes the proof.
\end{pf} \qed

By the definition of $H_t$ in (2.4), we have
\[\Im H_t(z)=(\Im
z)\left(1-(t-1)\int_\mathbb{R}\frac{s^2+1}{|z-s|^2}\;d\rho(s)\right),\;\;\;z\in\mathbb{C}^+.
\leqno{(3.1)}\]

\begin{prop} The function $H_t$ satisfies $\Im H_t(z)>0$ for all
$z\in\mathbb{C}^+$ if and only if $\mu$ is a Dirac measure in which
case $\mu^{\boxplus t}$ is a Dirac measure as well.
\end{prop}

\begin{pf} It follows from (3.1) that $\Im H_t(z)>0$ holds for all $z\in\mathbb{C}^+$
if and only if
\[\int_\mathbb{R}\frac{s^2+1}{|z-s|^2}\;d\rho(s)<\frac{1}{t-1},\;\;\;z\in\mathbb{C}^+,\]
which happens if and only if $g(x)\leq(t-1)^{-1}$ holds for all
$x\in\mathbb{R}$ by monotone convergence theorem. Then Lemma 3.1
shows that $\Im H_t(z)>0$, $z\in\mathbb{C}^+$, if and only if $\rho$
is a zero measure or, equivalently, $F_\mu=a+z$, i.e., $\mu$ is the
point mass $\delta_{-a}$. It is easy to see from (2.4) and (2.5)
that $F_{\mu^{\boxplus t}}(z)=ta+z$, whence $\mu^{\boxplus
t}=\delta_{-ta}$.
\end{pf} \qed

From the proof of the preceding proposition, $\mu$ is a point mass
if and only if $V_t^+$ is an empty set for all $t>1$. For the rest
of the paper, we confine our attention to the case of $\mu$ which is
not a Dirac measure.

The following lemmas are essentially similar to the results in [8]
and [11]. For the completeness, we provide the proofs here.

\begin{lem} The set $V_t^+$ coincides with
$\{x\in\mathbb{R}:f_t(x)>0\}$ and
$\{x\in\mathbb{R}:f_t(x)=0\}=V_t\cup V_t^-$. Moreover, for any
$x\in\mathbb{R}$ we have
\[\int_\mathbb{R}\frac{s^2+1}{(x-s)^2+f_t^2(x)}\;d\rho(s)\leq\frac{1}{t-1},\leqno{(3.2)}\] where the equality
holds for $x\in V_t^+$.
\end{lem}

\begin{pf} First note that if $z=x+iy\in\mathbb{C}^+$, i.e., $y>0$,
then
\[\Im
F_\mu(z)=y\left(1+\int_\mathbb{R}\frac{s^2+1}{(x-s)^2+y^2}\;d\rho(s)\right),\]
which gives, particularly,
\[\int_\mathbb{R}\frac{s^2+1}{(x-s)^2+y^2}\;d\rho(s)<\infty.\]
It follows that for any fixed $x\in\mathbb{R}$, the mapping
\[y\mapsto\int_\mathbb{R}\frac{s^2+1}{(x-s)^2+y^2}\;d\rho(s)\] is
decreasing and continuous on $(0,\infty)$ by Fatou's lemma. Indeed,
for any sequence $\{y_n\}$ converging to $y\in(0,\infty)$,
\begin{align*}
\int_\mathbb{R}\frac{s^2+1}{(x-s)^2+y^2}\;d\rho(s)&\leq
\liminf_{n\to\infty}\int_\mathbb{R}\frac{s^2+1}{(x-t)^2+y_n^2}\;d\rho(s) \\
&\leq\limsup_{n\to\infty}\int_\mathbb{R}\frac{s^2+1}{(x-s)^2+y_n^2}\;d\rho(s)\leq
\int_\mathbb{R}\frac{s^2+1}{(x-t)^2+y^2}\;d\rho(s)
\end{align*} which shows the continuity of the mapping.
If $f_t(x)>0$ then the definition of $f_t$ shows that for small
$\epsilon>0$ we have
\[\int_\mathbb{R}\frac{s^2+1}{(x-s)^2+(1+\epsilon)f_t^2(x)}\;d\rho(s)\leq
\frac{1}{t-1}<\int_\mathbb{R}\frac{s^2+1}{(x-s)^2+(1-\epsilon)f_t^2(x)}\;d\rho(s),\]
which gives the equality in (3.2) by letting $\epsilon\to0$. If
$f_t(x)=0$ then it is easy to see from Fatou's lemma that $x\in
V_t\cup V_t^-=\mathbb{R}\setminus V_t^+$, whence
$V_t^+=\{x\in\mathbb{R}:f_t(x)>0\}$ and $V_t\cup
V_t^-=\{x\in\mathbb{R}:f_t(x)=0\}$. This completes the proof.
\end{pf} \qed

With the help of preceding lemma, we are able to provide another
description of the set $\Omega_t=\{z\in\mathbb{C}^+:\Im H_t(z)>0\}$.

\begin{lem} The set $\{x+iy\in\mathbb{C}^+:y>f_t(x)\}$ coincides with the set $\Omega_t$ and
\[\int_\mathbb{R}\frac{s^2+1}{|z-s|^2}\;d\rho(s)<\frac{1}{t-1},\;\;\;\;z\in\Omega_t.\leqno{(3.3)}\]
The function $f_t$ is continuous on $\mathbb{R}$ and the set $V_t^+$
is open. Moreover, if $1<t_1<t_2$ then
$\Omega_{t_2}\subset\Omega_{t_1}$ and $f_{t_1}(x)\leq f_{t_2}(x)$
for all $x\in\mathbb{R}$.
\end{lem}

\begin{pf} The equation (3.1) shows that $z\in\Omega_t$ if and only if (3.3) holds.
By Lemma 3.3 and the definition of $f_t$, we have
$\Omega_t=\{x+iy\in\mathbb{C}^+:y>f_t(x)\}$, and hence the boundary
$\partial\Omega_t$ is the graph of $y=f_t(x)$. Since
$\partial\Omega_t$ is simply connected and
$V_t=\{x\in\mathbb{R}:f_t(x)>0\}$, the second assertion holds. The
last assertion follows from (3.3).
\end{pf} \qed

It was shown in [8] that the Cauchy transform $G_\mu$ is Lipschitz
continuous on some subset of $\overline{\mathbb{C}^+}$. The
following proposition gives similar results for the reciprocal
Cauchy transform $F_\mu$ on $\overline{\Omega_t}$.

\begin{prop} The reciprocal Cauchy transform $F_\mu$ extends continuously to
$\overline{\Omega_t}$ and satisfies
\[|F_\mu(z_1)-F_\mu(z_2)|\leq\frac{t}{t-1}|z_1-z_2|,\;\;\;\;z_1,z_2\in\overline{\Omega_t}.\leqno{(3.4)}\]
Moreover, this continuous extension can be represented as
\[F_\mu(z)=a+z+\int_\mathbb{R}\frac{1+sz}{s-z}\;d\rho(s),\;\;\;\;z\in\overline{\Omega_t}.\leqno{(3.5)}\] If
$t>2$ then
\[\frac{t-2}{t-1}|z_1-z_2|\leq|F_\mu(z_1)-F_\mu(z_2)|,\;\;\;z_1,z_2\in\overline{\Omega_t},\]
and, consequently, $F_\mu$ is one-to-one.
\end{prop}

\begin{pf} Using (3.3) and monotone convergence theorem gives the inequality
\[\int_\mathbb{R}\frac{s^2+1}{|z-s|^2}\;d\rho(s)\leq\frac{1}{t-1},\;\;\;\;z\in\overline{\Omega_t},\]
whence the integral
\[\int_\mathbb{R}\frac{s}{|z-s|^2}\;d\rho(s)\] converges
for all $z\in\overline{\Omega_t}$ by H\"{o}lder inequality. This
implies that
\[\int_\mathbb{R}\frac{1+sz}{s-z}\;d\rho(s)=\int_\mathbb{R}
\frac{(1-|z|^2)s+x(s^2-1)}{|s-z|^2}\;d\rho(s)+iy\int_\mathbb{R}
\frac{s^2+1}{|s-z|^2}\;d\rho(s)\] converges as well for
$z=x+iy\in\overline{\Omega_t}$. This implies that (3.5) holds, and
moreover for $z_1,z_2\in\overline{\Omega_t}$ we have
\begin{align*}
\left|\frac{F_{\mu}(z_1)-F_{\mu}(z_2)}{z_1-z_2}\right|&
\leq1+\left|\int_\mathbb{R}\frac{s^2+1}{(s-z_1)(s-z_2)}\;d\rho(s)\right| \\
&\leq1+\left(\int_\mathbb{R}\frac{s^2+1}{|s-z_1|^2}\;d\rho(s)\right)^{1/2}
\left(\int_\mathbb{R}\frac{s^2+1}{|s-z_2|^2}d\rho(s)\right)^{1/2} \\
&\leq1+\frac{1}{t-1}=\frac{t}{t-1},
\end{align*} where the H\"{o}lder inequality is used in the second
inequality. Similarly, if $t>2$ then
\[\left|\frac{F_{\mu}(z_1)-F_{\mu}(z_2)}{z_1-z_2}\right|\geq1-\frac{1}{t-1},\]
which implies the desired inequality. This completes the proof.
\end{pf} \qed

For $t>1$, it could happen that the boundary $\partial\Omega_t$
contains an interval $I\subset\mathbb{R}$, i.e., $f_t$ vanishes on
$I$. The following corollary characterizes such an interval.

\begin{cor} Let $I$ be an open interval in $\mathbb{R}$. Then
$\rho(I)=0$ if and only if for any closed interval $J\subset I$
there exists some $t>1$ such that $f_t(x)=0$ for $x\in J$. Moreover,
if $\rho(I)=0$ then the expression $(3.5)$ for $F_\mu$ holds for
$x\in I$ and $F_\mu$ is strictly increasing on $I$.
\end{cor}

\begin{pf} First, assume that
$I$ is bounded, $\rho(I)=0$, and $J$ is any closed interval
contained in $I$. Since $g$ is continuous on $I$, there exists a
number $t>1$ such that $g(x)\leq(t-1)^{-1}$ or, equivalently,
$f_t(x)=0$ for all $x\in J$. If $I=(a,\infty)$ (resp.
$I=(-\infty,a)$) for some finite $a$ and $\rho(I)=0$ then $g$ is
decreasing (resp. increasing) on $I$, whence it is easy to see that
the necessity follows. The sufficiency follows from Lemma 3.1 and
3.3.

By Proposition 3.5 and the first assertion of the corollary, it is
easy to see that (3.5) holds on any interval $I$ having
$\rho$-measure zero. For such an interval $I$, taking the derivative
gives that
\[F_\mu'(x)=1+\int_\mathbb{R}\frac{s^2+1}{(s-x)^2}\;d\rho(s),\;\;\;\;x\in I,\] and therefore
$F_\mu$ is strictly increasing on $I$.
\end{pf} \qed

From the proof of Proposition 3.5, the function $H_t$ also satisfies
similar estimates. That is,
\[|H_t(z_1)-H_t(z_2)|\leq2|z_1-z_2|,\;\;\;\;z_1,z_2\in\overline{\Omega_t},\]
from which we deduce that $H_t$ is also Lipscitz continuous and has
a continuous extension to $\overline{\Omega_t}$. If this continuous
extension is still denoted by $H_t$ then the equation
$H_t(\omega_t(z))=z$, $z\in\overline{\mathbb{C}^+}$, shows that
\[2|z_1-z_2|\leq|\omega_t(z_1)-\omega_t(z_2)|,\;\;\;z_1,z_2\in\overline{\mathbb{C}^+}.\]
The results in the following corollary are direct consequences of
(3.3) and preceding discussions. We refer the reader to [2] for
related results of the function $H_t$ and the subordination function
$\omega_t$.

\begin{cor} The function $H_t$ is a conformal mapping
from $\Omega_t$ onto $\mathbb{C}^+$ and is a homeomorphism from
$\overline{\Omega_t}$ to $\overline{\mathbb{C}^+}$.
\end{cor}

Note that the continuity of $f_t$ on $\mathbb{R}$ ensures that the
map $x\mapsto x+if_t(x)$ is a homeomorphism of $\mathbb{R}$ onto
$\partial\Omega_t$. Define the map $\psi_t:\mathbb{R}\to\mathbb{R}$
as
\[\psi_t(x)=H_t(x+if_t(x)).\] Then
$\psi_t$ is a homeomorphism as well by Corollary 3.7 and the
equation
\[\omega_t(\psi_t(x))=\omega_t(H_t(x+if_t(x)))=x+if_t(x)\;\;\;\mathrm{holds}\;\;\;\mathrm{for}
\;\;\;x\in\mathbb{R}.\] It is apparent that we have
\[F_{\mu^{\boxplus t}}(\psi_t(x))=\frac{tx-\psi_t(x)+itf_t(x)}{t-1},\;\;\;x\in\mathbb{R}.\leqno{(3.6)}\]
Since $\omega_t$ can be continued
analytically to a neighborhood of $\psi_t(x)$ if $f_t(x)>0$, it
follows that the function $f_t$ is analytic on $V_t^+$.

Now we are in a position to state the main theorem of this section.
For any measure $\nu$ on $\mathbb{R}$, let $\mathrm{supp}(\nu)$ be
the support of $\nu$ and denote by $\nu^{\mathrm{ac}}$ the
absolutely continuous part of $\nu$ with respect to Lebesgue
measure.

\begin{thm} Suppose that $\mu$ is a Borel probability measure $($not a point mass$)$ on $\mathbb{R}$ and
$t>1$. Then the following statements hold.
\\ $(i)$ The absolutely continuous part of $\mu^{\boxplus t}$ is concentrated on
the closure of $\psi_t(V_t^+)$.
\\ $(ii)$ The density of $(\mu^{\boxplus t})^{\mathrm{ac}}$ on the set $\psi_t(V_t^+)$ is given by
\[\frac{d\mu^{\boxplus t}}{dx}(\psi_t(x))=\frac{t(t-1)f_t(x)}{\pi|tx-\psi_t(x)+itf_t(x)|^2},\;\;\;x\in V_t^+.\]
$(iii)$ The density of $(\mu^{\boxplus t})^{\mathrm{ac}}$ is
analytic on the set $\psi_t(V_t^+)$ and at the point $\psi_t(x)$,
$x\in V_t^+$, it is bounded by $(t-1)/(\pi tf_t(x))$.
\\ $(iv)$ The number of the
components in $\mathrm{supp}((\mu^{\boxplus t})^{\mathrm{ac}})$ is a
decreasing function of $t$.
\end{thm}

\begin{pf} Since $F_{\mu^{\boxplus t}}$ is continuous on
$\overline{\mathbb{C}^+}$, the density of the nonatomic part of
$\mu^{\boxplus t}$ is continuous except at the points $x$ at which
$F_{\mu^{\boxplus t}}(x)=0$. In view of (3.6) and the Stieltjes
inversion formula, the density of $\mu^{\boxplus t}$ at $\psi_t(x)$
where $x\in V_t^+$ is given by
\begin{align*}
\frac{d\mu^{\boxplus t}}{dx}(\psi_t(x))&=-\frac{1}{\pi}\Im
G_{\mu^{\boxplus
t}}(\psi_t(x)) \\
&=\frac{t(t-1)f_t(x)}{\pi|tx-\psi_t(x)+itf_t(x)|^2}.
\end{align*}
This shows that the support of $(\mu^{\boxplus t})^{\mathrm{ac}}$ is
the closure of $\psi_t(V_t^+)$ and they have the same number of
components since $\psi_t$ is a homeomorphism. The assertion (iii)
follows from part (ii) and the fact that $f_t$ is analytic on
$V_t^+$. To verify the statement (iv), we need to show the the
number of the components in $V_t^+$ is nonincreasing as $t$
increases. This will hold if we show that $g$ never has a local
maximum in any open interval $(a,b)$ in $\mathbb{R}\backslash
V_t^+$. Indeed, $g$ is strictly convex on such an interval by Lemma
3.1, whence (iv) holds.
\end{pf} \qed

Before we state the next result, observe that the number of atoms of
$\mu^{\boxplus t}$ decreases as a function of $t$. Indeed, this is a
direct consequence of the fact that a point $\alpha$ is an atom of
$\mu^{\boxplus t}$ if and only $\mu(\{\alpha/t\})>1-t^{-1}$. Since
the open set $V_t^+$ can be written as a countable union of open
intervals, combining the result for atoms and Theorem 3.8 gives the
following conclusion.

\begin{cor} With the same assumption in Theorem $\mathrm{3.8}$, the measure $\mu^{\boxplus
t}$ has at most countably many components in the support, which
consists of finitely many points $($atoms$)$ and countably many
intervals. Moreover, the number of the components in
$\mathrm{supp}(\mu^{\boxplus t})$ is a decreasing function of $t$.
\end{cor}

\section{Support and Regularity for $\mu^{\boxplus t}$}

In this section, we will investigate the support and regularity
property for the free convolution power $\mu^{\boxplus t}$ of $\mu$
where $t>1$ and $\mu$ is a Borel probability measure (not a point
mass so that $V_t^+$ is nonempty) on $\mathbb{R}$. By analyzing the
formula
\[F_{\mu^{\boxplus t}}(\psi_t(x))=\frac{tx-\psi_t(x)+itf_t(x)}{t-1},\;\;\;x\in\mathbb{R},\leqno{(4.1)}\]
established in the previous section, we are able to understand the
relation between supports and regularities for the measures $\mu$
and $\mu^{\boxplus t}$. If $f_t(x)=0$ for some $x\in\mathbb{R}$ then
$H_t(x)=\psi_t(x)$ and (4.1) can be rewritten as
\[F_\mu(x)=F_{\mu^{\boxplus
t}}(\psi_t(x))=\frac{tx-\psi_t(x)}{t-1}.\leqno{(4.2)}\] The
following result gives another basic property of the homeomorphism
$\psi_t$.

\begin{lem} The derivative $\psi_t'(x)>0$
for $x\in V_t^+$ and $\psi_t$ is a strictly increasing function on
$\mathbb{R}$.
\end{lem}

\begin{pf} Let $z=x+if_t(x)$ be any point with $x\in V_t^+$. Since functions $H_t$ and $f_t$ are
analytic at points $z$ and $x$, respectively, it follows that
\begin{align*}
\psi_t'(x)&=\left(\frac{dH_t}{dz}(x+if_t(x))\right)(1+if_t'(x)) \\
&=\Re\left(\frac{dH_t}{dz}(x+if_t(x))\right)-\Im\left(\frac{dH_t}{dz}(x+if_t(x))\right)f_t'(x).
\end{align*} Observe that
\begin{align*}
\frac{dH_t}{dz}(x+if_t(x))&=1-(t-1)\int_\mathbb{R}\frac{s^2+1}{(s-z)^2}\;d\rho(s) \\
&=1-(t-1)\int_\mathbb{R}\frac{(s^2+1)[(s-x)^2-f_t^2(x)+2i(s-x)f_t(x)]}{|s-z|^4}\;d\rho(s) \\
&=1-(t-1)\int_\mathbb{R}\frac{(s^2+1)[|s-z|^2-2f_t^2(x)+2i(s-x)f_t(x)]}{|s-z|^4}\;d\rho(s) \\
&=2(t-1)\left[\int_\mathbb{R}\frac{(s^2+1)f_t^2(x)}{|s-z|^4}\;d\rho(s)-i
\int_\mathbb{R}\frac{(s^2+1)(s-x)f_t(x)}{|s-z|^4}\;d\rho(s)\right]
\end{align*} where the identity
\[\int_\mathbb{R}\frac{s^2+1}{|s-z|^2}\;d\rho(s)=\frac{1}{t-1}\leqno{(4.3)}\] is used in
the forth equality. On the other hand, differentiating (4.3) with
respect to $x$ gives
\[\int_\mathbb{R}\frac{(s^2+1)(s-x)}{|s-z|^4}\;d\rho(s)=f_t(x)f_t'(x)
\int_\mathbb{R}\frac{s^2+1}{|s-z|^4}\;d\rho(s),\] whence
\[\frac{dH_t}{dz}(x+if_t(x))=2(t-1)f_t^2(x)\left[\int_\mathbb{R}\frac{s^2+1}{|s-z|^4}\;d\rho(s)-if_t'(x)
\int_\mathbb{R}\frac{s^2+1}{|s-z|^4}\;d\rho(s)\right].\] Combining
these results gives that
\begin{align*}
\psi_t'(x)&=2(t-1)f_t^2(x)\left[\int_\mathbb{R}\frac{s^2+1}{|s-z|^4}\;d\rho(s)+f_t'^2(x)
\int_\mathbb{R}\frac{s^2+1}{|s-z|^4}\;d\rho(s)\right] \\
&=2(t-1)f_t^2(x)(f_t'^2(x)+1)\int_\mathbb{R}\frac{s^2+1}{|s-z|^4}\;d\rho(s) \\
&>0,
\end{align*} as desired. Since $\psi_t$ is one-to-one and
continuous on $\mathbb{R}$ with nonnegative derivative on a nonempty
subset of $\mathbb{R}$, it must be strictly increasing.
\end{pf} \qed

Note that if $x\in V_t$ then $f_t(x)=0$, and Proposition 3.5 shows
that the limit $F_\mu(x)=\lim_{\epsilon\to0}F_\mu(x+i\epsilon)$
exists and is a finite number. The next result characterizes the
atoms of $\mu$ in terms of $V_t$ and $F_\mu$.

\begin{lem} A point $\alpha\in\mathbb{R}$ is an
atom of $\mu$ with mass $\mu(\{\alpha\})=1-t^{-1}$ if and only if
$\alpha\in V_t$ and $F_\mu(\alpha)=0$ in which case the
Julia-Carath\'{e}odory derivative of $F_\mu$ at $\alpha$ is
\[F_\mu'(\alpha)=\frac{1}{\mu(\{\alpha\})}=1+\int_\mathbb{R}\frac{s^2+1}{(s-\alpha)^2}\;d\rho(s).\]
\end{lem}

\begin{pf} First assume that $\alpha$ is an atom of $\mu$, i.e., the
limit in (2.3) is $1-t^{-1}>0$. This implies that
\begin{align*}
\lim_{\epsilon\to0}\frac{\Re F_\mu(\alpha+i\epsilon)}{i\epsilon}&
=\frac{1}{2}\lim_{\epsilon\to0}\left(\frac{F_\mu(\alpha+i\epsilon)}{i\epsilon}
+\frac{\overline{F_\mu(\alpha+i\epsilon)}}{i\epsilon}\right) \\
&=\frac{1}{2}\lim_{\epsilon\to0}\left(\frac{1}{i\epsilon
G_\mu(\alpha+i\epsilon)}-\frac{1}{\overline{i\epsilon
G_\mu(\alpha+i\epsilon)}}\right) \\
&=0.
\end{align*}
Using the limit shown above gives
\begin{align*}
\frac{1}{\mu(\{\alpha\})}&=\lim_{\epsilon\to0}\frac{F_\mu(\alpha+i\epsilon)}{i\epsilon} \\
&=\lim_{\epsilon\to0}\frac{\Im F_\mu(\alpha+i\epsilon)}{\epsilon} \\
&=\lim_{\epsilon\to0}\frac{\epsilon+\epsilon\int_\mathbb{R}\frac{s^2+1}{(s-\alpha)^2+\epsilon^2}\;d\rho(s)}{\epsilon} \\
&=1+\lim_{\epsilon\to0}\int_\mathbb{R}\frac{s^2+1}{(s-\alpha)^2+\epsilon^2}\;d\rho(s) \\
&=1+\int_\mathbb{R}\frac{s^2+1}{(s-\alpha)^2}\;d\rho(s),
\end{align*} where in the last equality monotone convergence theorem is used.
Hence $\alpha$ belongs to the set $V_t$, and $F_\mu(\alpha)$ is
defined and equals zero.

To show the sufficiency, by the arguments shown above it suffices to
show that
\[\lim_{\epsilon\to0}\frac{\Re F_\mu(\alpha+i\epsilon)}{i\epsilon}=0\]
or, equivalently,
\[\lim_{\epsilon\to0}\frac{\Re[F_\mu(\alpha+i\epsilon)-F_\mu(\alpha)]}{\epsilon}=0.\leqno{(4.4)}\]
First, in view of (3.5) we obtain
\[\frac{F_\mu(\alpha+i\epsilon)-F_\mu(\alpha)}{\epsilon}=i\left(1+\int_\mathbb{R}\frac{s^2+1}{(s-\alpha-i\epsilon)
(s-\alpha)}\;d\rho(s)\right).\] Since
\[\int_\mathbb{R}\left|\frac{s^2+1}{(s-\alpha-i\epsilon)
(s-\alpha)}\right|d\rho(s)\leq\int_\mathbb{R}\frac{s^2+1}{(s-\alpha)^2}\;d\rho(s)=\frac{1}{t-1},\]
it follows that
\[\lim_{\epsilon\to0}\frac{F_\mu(\alpha+i\epsilon)-F_\mu(\alpha)}{\epsilon}=i\frac{t}{t-1}\]
by dominated convergence theorem which implies (4.4). The expression
of the Julia-Carath\'{e}odory derivative of $F_\mu$ is clear.
\end{pf} \qed

It was proved in [2] that if a point $x\in\mathbb{R}$ belongs to
$\omega_t(\mathbb{R})$, i.e., $f_t(x)=0$ then the
Julia-Carath\'{e}odory derivative $H_t'(x)\geq0$, where the
existence of
\[H_t'(x)=\lim_{\epsilon\downarrow0}\frac{H(x+i\epsilon)-H(x)}{i\epsilon}\]
is guaranteed by the existence of $H_t(x)\in\mathbb{R}$. Indeed, by
a similar proof in Lemma 4.2 and Lemma 3.3 at such a point $x$ the
Julia-Carath\'{e}odory derivative of $F_\mu$ is given by
\[F_\mu'(x)=1+\int_\mathbb{R}\frac{s^2+1}{(s-x)^2}\;d\rho(s)\leq
t/(t-1),\leqno{(4.5)}\] from which we deduce that
\[H_t'(x)=1-(t-1)\int_\mathbb{R}\frac{s^2+1}{(s-x)^2}\;d\rho(s)\geq0.\leqno{(4.6)}\]
Note that (4.6) shows that if $f_t(x)=0$ then $H_t'(x)>0$ if and
only if $x\in V_t^-$. In the following proposition, by means of
functions $f_t$ and $\psi_t$ we are able to obtain certain
regularity results which were proved in [2].

\begin{prop} Let $\alpha$ be a point in $\mathbb{R}$ and $t>1$. Then the following conditions are
equivalent:
\\ $(i)$ $f_t(\alpha)=0$ and $F_\mu(\alpha)=0$;
\\ $(ii)$ $F_{\mu^{\boxplus t}}(t\alpha)=0$;
\\ $(iii)$ $\mu(\{\alpha\})\geq1-t^{-1}$.
\end{prop}

\begin{pf} The equivalence of (i) and (ii) follows from the bijectivity of
$\psi_t$ and (4.1). Indeed, letting $\psi_t(\alpha')=t\alpha$ for
some $\alpha'\in\mathbb{R}$ implies that $F_{\mu^{\boxplus
t}}(t\alpha)=0$ if and only if $t\alpha'=\psi_t(\alpha')$ and
$f_t(\alpha')=0$ which translate to the conditions in (i). By Lemma
4.2 and (4.5), it is easy to see that (i) implies (iii). On the
other hand, by Lemma 4.2 we see that if
$\mu(\{\alpha\})=1-t_0^{-1}\geq1-t^{-1}$ then $F_\mu(\alpha)=0$ and
$f_{t_0}(\alpha)=0$. Since $t\leq t_0$, $f_t(\alpha)=0$ by Lemma
3.4, whence (iii) implies (i).  This completes the proof.
\end{pf} \qed

Next, we investigate those intervals on which $f_t$ vanishes.

\begin{prop} If $f_t$ vanishes on some interval $I$ then
the following statements are equivalent:
\\ $(i)$ $\mu(I)=0$;
\\ $(ii)$ $F_\mu(x)\neq0$ for any $x$ in $I$.
\\ If $(i)$ or $(ii)$ is satisfied then $\mu^{\boxplus
t}(\psi_t(I))=0$.
\end{prop}

\begin{pf}  The implication that (i) implies (ii) follows from Proposition 4.3.
Indeed, if $\Re F_\mu(x)=0$ for some $x\in I$ then $x$ is an atom of
$\mu$, and hence $\mu(I)>0$. Next, we prove that (ii) implies (i).
First we consider the case that $I$ is a compact set $[a,b]$. Since
$F_\mu$ is continuous and nonzero on $[a,b]$, $G_\mu(x)$ is a
continuous and real-valued function on $[a,b]$, and consequently
$\mu([a,b])=0$ by the Stieltjes inversion formula. Next, suppose
that $I$ is bounded but not closed. If $I=(a,b]$ then we have
$\mu([a+\epsilon,b])=0$ for any small $\epsilon>0$ by the result
established above, and therefore letting $\epsilon\to0$ gives
$\mu((a,b])=0$. Similarly, the assertion (i) holds if $I=[a,b)$ or
$I=(a,b)$. Finally, if $I$ is unbounded (closed or not closed) then
by the results shown above and the countable additivity of $\mu$, it
is easy to see that $\mu(I)=0$. This completes the proof that (ii)
implies (i).

By similar arguments and (4.2), it is easy to see that the last
statement follows from (i) or (ii).
\end{pf} \qed

\begin{cor} Assume that $I\subset\mathbb{R}$ is an interval on which $f_t$ vanishes and $\mu(I)>0$.
Then the interval $I$ contains one and only one atom $\alpha$ of
$\mu$ and this only atom has mass $\mu(\{\alpha\})\geq1-t^{-1}$.
Moreover, $\mu(I\setminus\{\alpha\})=0$ and $\mu^{\boxplus
t}(\psi_t(I))=t\mu(\{\alpha\})-(t-1)$.
\end{cor}

\begin{pf} First observe that $F_\mu$ has zeros in $I$. Since $F_\mu$
is strictly increasing on $I$ by Corollary 3.6, it only has one zero
$\alpha$, i.e., $\alpha$ is the only atom of $\mu$ in $I$ by Lemma
4.2. The atom $\alpha$ with the desired mass is a direct consequence
of Proposition 4.3. It is clear that the set $I\backslash\{\alpha\}$
must have $\mu$-measure zero. Since $I\backslash\{\alpha\}\subset
V_t\cup V_t^-$ and $\mu^{\boxplus
t}(\{t\alpha\})=t\mu(\{\alpha\})-(t-1)$, $\psi_t(I)$ has the desired
$\mu^{\boxplus t}$-measure.
\end{pf} \qed

\begin{thm} Let $I$ be any component in $\mathrm{supp}(\mu)$ and $t>1$. If $I$ is an
interval then it does not contain a closed interval on which $f_t$
vanishes and the interval $\psi_t(I)$ is contained in some component
of $\mathrm{supp}(\mu^{\boxplus t})$.
\end{thm}

\begin{pf} Let $\alpha_1,\cdots,\alpha_n$ be all the atoms
of $\mu$ in $I$ with $\mu(\{\alpha_j\})\geq1-t^{-1}$,
$j=1,\cdots,n$. Then it suffices to show that there does not exist
an open interval $J\subset I\backslash\{\alpha_1,\cdots,\alpha_n\}$
such that $f_t$ vanishes on $J$. If such an interval $J$ exists then
the fact $\mu(J)>0$ gives the existence of an atom $\alpha$ of $\mu$
in $J$ with mass $\mu(\{\alpha\})\geq1-t^{-1}$ by Corollary 4.5, a
contradiction. Hence the first desired result follows. This result
also reveals that $g>(t-1)^{-1}$ on the interval $I$ except at
points $x$ such that $f_t(x)=0$. Since $\mu^{\boxplus t}$ does not
contain the singular continuous part, the last assertion follows
from Theorem 3.8(i).
\end{pf} \qed

It is known that the measure $\mu^{\boxplus t}$ has fewer atoms when
$t$ increases. The next result explains where these disappearing
atoms go.

\begin{prop} Assume that $\alpha$ is an atom of $\mu$ and
$\mu(\{\alpha\})=1-t_0^{-1}$. Then $\mu^{\boxplus
t}(\{\psi_t(\alpha)\})>0$ for $1<t<t_0$ and
$\psi_t(\alpha)\in\mathrm{supp}((\mu^{\boxplus t})^{\mathrm{ac}})$
for $t\geq t_0$.
\end{prop}

\begin{pf} For $1<t<t_0$, we have $\mu(\{\alpha\})>1-t^{-1}$, which implies
that $\psi_t(\alpha)=t\alpha$ is an atom of $\mu^{\boxplus t}$ by
(4.1) and Proposition 4.3. Next note that there does not exist an
open interval $I$ containing $\alpha$ such that $f_{t_0}(x)=0$ for
all $x\in I$. Indeed, if such an interval $I$ exists then the facts
that $g(\alpha)=(t_0-1)^{-1}$ and $g$ is strictly convex on $I$ will
lead to a contradiction. Hence $\psi_{t_0}(\alpha)$ must be in the
closure of $\psi_{t_0}(V_{t_0}^+)$. Since $\alpha\in V_t^+$ for
$t>t_0$, the desired result follows from Theorem 3.8(i).
\end{pf} \qed

Theorem 4.6 and Proposition 4.7 give the inclusion
$\psi_t(\mathrm{supp}(\mu))\subset\mathrm{supp}(\mu^{\boxplus t})$
or, equivalently, the complement
$\mathbb{R}\backslash\mathrm{supp}(\mu^{\boxplus t})$ is contained
in $\psi_t(\mathbb{R}\backslash\mathrm{supp}(\mu))$ for all $t>1$.
Since $\mathbb{R}\backslash\mathrm{supp}(\mu)$ is a countable union
of open intervals, the preceding observation leads us to investigate
open intervals $I$ which have $\mu$-measure zero. First note that
the Cauchy transform $G_\mu$ extends analytically through the
interval $I$ and takes real values on $I$. Moreover, it can be
expressed as
\[G_\mu(x)=\int_\mathbb{R}\frac{1}{x-s}\;d\mu(s),\;\;\;x\in I,\]
and $G_\mu'(x)=-\int(s-x)^{-2}d\mu(s)<0$ which shows that $G_\mu$
has at most one zero in $I$. To analyze the values of $f_t$ on the
set $\mathbb{R}\backslash\mathrm{supp}(\mu)$, we need the following
lemma.

\begin{lem} Suppose that $I$ is an open interval in $\mathbb{R}$ with $\mu(I)=0$.
\\ $(i)$ If $G_\mu$ has a zero $x_0$ in $I$ then $x_0$ is an atom of
$\rho$ with mass
\[\rho(\{x_0\})=\frac{-1}{(x_0^2+1)G_\mu'(x_0)}\] and
$\rho(I\setminus\{x_0\})=0$.
\\ $(ii)$ At any point $x\in I$, the function $g$ can be expressed as
\[g(x)=\frac{-G_\mu'(x)}{G_\mu^2(x)}-1,\] where the
equality is interpreted as $+\infty$ on both sides at the unique
zero $x_0$, if it exists, of $G_\mu$ on $I$. In addition, if
$I=(a,b)$ is bounded then
\[\inf_{x\in I}g(x)\geq\frac{\rho(\mathbb{R})}{1+\max\{a^2,b^2\}}.\]
\end{lem}

\begin{pf} First observe that $F_\mu$ is continuous and takes real values on
$I\backslash\{x_0\}$, whence the weak$^*$-limit of the measures
$d\rho_\epsilon(s)$ in (2.2) as $\epsilon\downarrow0$ is zero which
gives $\rho(I\backslash\{x_0\})=0$. By Lemma 2.1, we have
\[\lim_{\epsilon\downarrow0}i\epsilon
F_\mu(x_0+i\epsilon)=-(x_0^2+1)\rho(\{x_0\}).\] On the other hand,
\[\lim_{\epsilon\downarrow0}
\frac{G_\mu(x_0+i\epsilon)}{i\epsilon}=\lim_{\epsilon\downarrow0}
\frac{G_\mu(x_0+i\epsilon)-G_\mu(x_0)}{i\epsilon}=G_\mu'(x_0),\]
which gives the desired mass for $\rho(\{x_0\})$ in (i). To verify
the statements in (ii), by (i) we may assume that $G_\mu$ is nonzero
on $I$. Then observe that for any $x\in I$ and $\epsilon>0$ we have
\begin{align*}
\Im
F_\mu(x+i\epsilon)&=\Im\frac{\overline{G_\mu(x+i\epsilon)}}{|G_\mu(x+i\epsilon)|^2}
=\frac{-\Im
G_\mu(x+i\epsilon)}{|G_\mu(x+i\epsilon)|^2}
=\epsilon\frac{\int_\mathbb{R}\frac{d\mu(s)}{(x-s)^2+\epsilon^2}}{\left|\int_\mathbb{R}\frac{d\mu(s)}
{x-s+i\epsilon}\right|^2}.
\end{align*} On the other hand, the Nevanlinna representation for
$F_\mu$ gives
\[\Im
F_\mu(x+i\epsilon)=\epsilon\left(1+\int_\mathbb{R}\frac{s^2+1}{(s-x)^2+\epsilon^2}\;d\rho(s)\right),\]
from which we obtain
\[\int_\mathbb{R}\frac{s^2+1}{(s-x)^2+\epsilon^2}\;d\rho(s)
=\frac{\int_\mathbb{R}\frac{d\mu(s)}{(x-s)^2+\epsilon^2}}{\left|\int_\mathbb{R}\frac{d\mu(s)}
{x-s+i\epsilon}\right|^2}-1.\leqno{(4.7)}\] Since $\mu(I)=0$ and
$G_\mu$ is nonzero on $I$, it follows that
\[\lim_{\epsilon\to0}\int_\mathbb{R}\frac{d\mu(s)}{(x-s)^2+\epsilon^2}=
\int_\mathbb{R}\frac{d\mu(s)}{(x-s)^2}<\infty\] and
\[\lim_{\epsilon\to0}\int_\mathbb{R}\frac{d\mu(s)}{(x-s)+i\epsilon}=
\int_\mathbb{R}\frac{d\mu(s)}{x-s}\neq0\] for all $x\in I$.
Therefore, by the equation (4.7) we have
\begin{align*}
g(x)&=\int_\mathbb{R}\frac{s^2+1}{(s-x)^2}\;d\rho(s) \\
&=\lim_{\epsilon\to0}\int_\mathbb{R}\frac{s^2+1}{(s-x)^2+\epsilon^2}\;d\rho(s) \\
&=\lim_{\epsilon\to0}\left(\frac{\int_\mathbb{R}\frac{d\mu(s)}{(x-s)^2+\epsilon^2}}{\left|\int_\mathbb{R}\frac{d\mu(s)}
{x-s+i\epsilon}\right|^2}-1\right) \\
&=\frac{\int_\mathbb{R}\frac{d\mu(s)}{(x-s)^2}}{\left(\int_\mathbb{R}\frac{d\mu(s)}
{x-s}\right)^2}-1 \\
&=\frac{-G_\mu'(x)}{G_\mu^2(x)}-1.
\end{align*} Finally, for any $x$ in the bounded interval $(a,b)$ and
$s\in\mathbb{R}\backslash(a,b)$ we have
\[\frac{s^2+1}{(s-x)^2}\geq\frac{1}{1+x^2}\geq\frac{1}{1+\max\{a^2,b^2\}},\]
from which we deduce that
\begin{align*}
g(x)&=\int_{-\infty}^a\frac{s^2+1}{(s-x)^2}\;d\rho(s)+
\int_b^\infty\frac{s^2+1}{(s-x)^2}\;d\rho(s) \\
&\geq\frac{\rho(\mathbb{R})}{1+\max\{a^2,b^2\}},
\end{align*} as desired. This finishes the proof.
\end{pf} \qed

\begin{prop} If $I$ is a bounded component in
$\mathbb{R}\backslash\mathrm{supp}(\mu)$ then for each $t>1$ the
interval $\psi_t(I)$ contains at most two components of
$\mathbb{R}\backslash\mathrm{supp}(\mu^{\boxplus t})$.
\end{prop}

\begin{pf} If $G_\mu$ has no zero in $I$ then the strict convexity of $g$ on $I$
shows that there exists at most one subinterval $J$ of $I$ such that
$f_t$ vanishes on $J$ when $t$ varies from one to infinity. If
$G_\mu$ has a zero $x_0$ in $I=(a,b)$ then applying the preceding
argument to the intervals $(a,x_0)$ and $(x_0,b)$ gives that two
subintervals at most. Hence the desired conclusion follows from
Theorem 3.8(i).
\end{pf} \qed

For any $t>1$, denote by $n(t)$ the number of components in
$\mathrm{supp}(\mu^{\boxplus t})$.

\begin{thm} If $\mu$ is any Borel probability measure on
$\mathbb{R}$ then the following statements are equivalent:
\\ $(i)$ $n(t)=1$ for sufficiently large $t$;
\\ $(ii)$ $n(t)<\infty$ for some $t>1$;
\\ $(iii)$ the infimum $m$ of $g$ on the set of all bounded
components of $\mathbb{R}\setminus\mathrm{supp}(\mu)$ is nonzero.
\\ Moreover, if $m>0$ then $n(t)=1$ for $t>t_0$, where
\[t_0=\max\left\{1+\frac{1}{m},\frac{1}{1-\mu(\{\alpha\})}\right\}\] and $\alpha$ is one of the atoms of
$\mu$ with the largest mass.
\end{thm}

\begin{pf} It is clear that (i) implies (ii). Next, suppose that
$n(t_1)<\infty$ for some $t_1>1$ and $\mu(\{\alpha\})<1-t_1^{-1}$
for any atom $\alpha$ of $\mu$. Let $S_{t_1}$ be the set of all
components which are bounded intervals in
$\{x\in\mathbb{R}:f_{t_1}(x)=0\}$. Then $S_{t_1}$ is a finite set
and by Corollary 4.5 the set $S_{t_1}$ is contained in the union of
some components $I_1,\cdots,I_n$ of
$\mathbb{R}\backslash\mathrm{supp}(\mu)$. Note that each $I_k$ is
bounded. Indeed, if $\mu((M,\infty))=0$ for some finite number $M$
then $G_\mu>0$ on $(M,\infty)$, whence $\rho((M,\infty))=0$ and $g$
is strictly decreasing on $(M,\infty)$. This implies that
$g(x)\leq(t_1-1)^{-1}$ for sufficiently large $x$, and so $f_{t_1}$
vanishes on $(c,\infty)$ for some finite number $c$. Similarly, no
$I_k$ is of the form $(-\infty,M)$. Let $[a,b]$ be the smallest
closed interval containing these $I_k$'s. Then
\[\inf\{g(x):x\in(\mathbb{R}\backslash\mathrm{supp}(\mu))\cap[a,b]\}\geq\frac{\rho(\mathbb{R})}{1+\max\{a^2,b^2\}}\]
by Lemma 4.8(ii). If $I$ is the intersection of any bounded
component in $\mathbb{R}\backslash\mathrm{supp}(\mu)$ with
$(b,\infty)$ and $I\neq\emptyset$ then $g(x)\geq(t_1-1)^{-1}$ on
$I$. Indeed, if $g(x)<(t_1-1)^{-1}$ for some $x\in I$ then there
exists some interval $J\subset I$ such that $g<(t_1-1)^{-1}$ on $J$,
i.e., $f_{t_1}=0$ on $J$, which violates the definition of
$S_{t_1}$. Similarly, $g\geq(t_1-1)^{-1}$ on any bounded component
in $(\mathbb{R}\backslash\mathrm{supp}(\mu))\cap(-\infty,a)$, and
hence (ii) implies (iii). Finally, assume that the number $m$ in the
statement (iii) is nonzero and let $t>t_0$, where $t_0$ is defined
as in theorem. Then with the help of Theorem 4.6, we see that
$\{x\in\mathbb{R}:f_{t}(x)=0\}$ contains at most two components
which are intervals, in which case these two intervals must be
unbounded. In other words, the closure of $V_t^+$ is an interval,
whence $n(t)=1$ and (iii) implies (i).
\end{pf} \qed

In the proof of the preceding theorem, the set
$\{x\in\mathbb{R}:f_t(x)=0\}$ might contain two unbounded components
for any $t>1$. Actually, this happens when $\mu$ is compactly
supported which was proved in [6] and [14]. In the following, we
provide an analytic way to prove this fact.

\begin{cor} Let $t>1$. Then the following statements hold.
\\ $(i)$ The function $f_t$ vanishes on $(a,\infty)$ $($resp.
$(-\infty,a)$$)$ for some finite number $a$ if and only if
$\mu((M,\infty))=0$ $($resp. $(-\infty,M)$$)$ for some finite number
$M$.
\\ $(ii)$ The measure $\mu$ has compact support if and only if so
does $\mu^{\boxplus t}$.
\\ $(iii)$ If $\mu$ is compactly supported then $n(t)=1$
for large $t$.
\end{cor}

\begin{pf} If $f_t(x)=0$ for $x\in(a,\infty)$ then by Corollary 4.5 there exists a point
$M\geq a$ such that $\mu((M,\infty))=0$. Similarly, if $f_t$
vanishes on $(-\infty,a)$ then $\mu((-\infty,M))=0$ for some finite
number $M\leq a$. Conversely, if $\mu((b,\infty))=0$ for some $b$
then as shown in the preceding theorem we have $f_t(x)=0$ for
sufficiently large $x$. Similarly, if $\mu((-\infty,b))=0$ then
$f_t(x)=0$, $x\leq a$, for some finite number $a$, and so (i)
follows. The assertion (ii) follows from (i), while (iii) is a
direct consequence of Lemma 4.8(ii) and Theorem 4.10.
\end{pf} \qed

Now we are able to give a complete description of
$\mathrm{supp}(\mu^{\boxplus t})$ in terms of $\mathrm{supp}(\mu)$.
Consider a component $I=(a,b)$ in
$\mathbb{R}\backslash\mathrm{supp}(\mu)$. The case $b=\infty$ is
discussed in Corollary 4.11. If $b<\infty$ then one of the following
cases holds.
\\ $(i)$ The set $\mathbb{R}\backslash\mathrm{supp}(\mu)$ contains a component $(b,c)$ in
which case $b$ is an atom.
\\ $(ii)$ There exists a point $c>b$ such that $(b,c)$ is contained in a
component of $\mathrm{supp}(\mu)$.
\\ $(iii)$ There exists a sequence of components $(a_n,b_n)$ in $\mathbb{R}\backslash\mathrm{supp}(\mu)$
such that $b_{n+1}\leq a_n$ for all $n$ and $a_n\downarrow b$ as
$n\to\infty$.
\\ Cases (i) and (ii) are
analyzed in Proposition 4.3, Theorem 4.6, and Proposition 4.9.
Applying Corollary 4.5 to case (iii), it is easy to see that for
each $t>1$ there exists a point $c>b$ such that there is no interval
contained in $(b,c)$ on which $f_t$ vanishes, whence
$g\geq(t-1)^{-1}$ on $(b,c)$ and $\psi_t((b,c))$ is contained in
some component of $\mathrm{supp}(\mu^{\boxplus t})$. Similar
statements also hold for $a$.

It was proved in [1] that $\mu^{\boxplus t}((t\alpha,t\beta))>0$ if
$t\alpha,t\beta$ are atoms of $\mu^{\boxplus t}$. Using functions
$f_t$ and $\psi_t$, we can prove a somewhat stronger result.

\begin{prop} If $\mu$ has atoms $\alpha<\beta$ such that
$\mu(\{\alpha\}),\mu(\{\beta\})\geq1-t^{-1}$ then $\mu^{\boxplus
t}((t\alpha,t\beta))>0$.
\end{prop}

\begin{pf} First observe that $f_t(\alpha)=f_t(\beta)=F_\mu(\alpha)=F_\mu(\beta)=0$
and $\psi_t(\alpha)=t\alpha$, $\psi_t(\beta)=t\beta$. Without loss
of generality, we assume that $F_{\mu^{\boxplus t}}(x)\neq0$ for all
$x$ in $(t\alpha,t\beta)$ or, equivalently, $F_{\mu^{\boxplus
t}}(\psi_t(x))\neq0$ for all $x\in(\alpha,\beta)$. By Theorem
3.8(i), it suffices to show that $f_t$ is positive at some point in
$[\alpha,\beta]$. If $f_t$ vanishes on $[\alpha,\beta]$ then $F_\mu$
is continuous and strictly increasing on $[\alpha,\beta]$ by
Proposition 3.5 and Corollary 3.6, a contradiction. This completes
the proof.
\end{pf} \qed

By the work of Maassen [13], for any probability measure $\nu$ on
$\mathbb{R}$ there exists a unique probability measure $\mu$ on
$\mathbb{R}$ with mean zero and unit variance such that
\[F_\mu(z)=z-G_\nu(z),\;\;\;z\in\mathbb{C}^+.\leqno{(4.8)}\] If $\rho$ denotes the finite positive
Borel measure in the Nevanlinna representation of $F_\mu$ as before
then
\[(s^2+1)\;d\rho(s)=d\nu(s)\] by the Stieltjes inversion formula and (2.2). In
this case,
\[H_t(z)=z+(t-1)G_\nu(z),\;\;\;\;z\in\mathbb{C}^+,\] and hence $z\in\Omega_t$ if and only if
\[\int_\mathbb{R}\frac{d\nu(s)}{|z-s|^2}<\frac{1}{t-1}.\] Moreover, the function
$g$ defined in section 3 is given by the formula
\[g(x)=\int_\mathbb{R}\frac{d\nu(s)}{(s-x)^2},\;\;\;x\in\mathbb{R}.\leqno{(4.9)}\]

\begin{prop} There exists a Borel probability measure $\mu$ on $\mathbb{R}$ such that the number of
components in $\mathrm{supp}(\mu^{\boxplus t})$ is infinite for any
$t>1$.
\end{prop}

\begin{pf} Define the probability measure $\nu$ on $\mathbb{R}$ as
\[\nu=\sum_{n=1}^\infty2^{-n}\delta_{2^n},\] and let $\mu$ be the
unique measure satisfying the requirement (4.8). Then (4.9) shows
that the function $g$ is continuous and strictly convex on the
intervals $I_n=(2^n,2^{n+1})$, $n\geq1$, and
\[\lim_{x\downarrow2^n}g(x)=\lim_{x\uparrow2^{n+1}}g(x)=\infty,\leqno{(4.10)}\] from which
we see that $g$ reaches its minimum $m_n$ on each $I_n$. If $x_n$ is
the middle point of $I_n$ then for each $s\in\mathbb{R}\backslash
I_n$,
\[(s-x_n)^2\geq\left(\frac{2^{n+1}-2^n}{2}\right)^2=2^{2n-2},\] from
which we deduce that
\[m_n\leq g(x_n)=\int_{\mathbb{R}\backslash
I_n}\frac{d\nu(s)}{(s-x_n)^2}\leq2^{-2n+2}.\] This shows that for
any $t>1$, there are infinitely many intervals
$\{I_{n_k}\}_{k=1}^\infty$ each of which contains a closed
subinterval $J_{n_k}$ such that $f_t$ vanishes on $J_{n_k}$. Note
that (4.10) shows that these intervals $J_{n_k}$ are separated by
intervals on which $f_t$ is positive. Hence we conclude that the
closure of the set
\[V_t^+=\{x\in\mathbb{R}:g(x)>(t-1)^{-1}\}\] contains infinitely
many components, whence the desired result follows.
\end{pf} \qed

We conclude this section with an example. For any $0<\epsilon<1$,
consider the measure
\[\mu_\epsilon=\frac{\epsilon}{2}(\delta_{-1}+\delta_1)+(1-\epsilon)\delta_0.\]
Since
\[G_{\mu_\epsilon}(z)=\frac{z^2-1+\epsilon}{z(z^2-1)}\] has zeros at
$\pm\sqrt{1-\epsilon}$, by Lemma 4.8 and some simple manipulations
we obtain
\[\rho_\epsilon=\frac{\epsilon}{2(2-\epsilon)}\left(\delta_{-\sqrt{1-\epsilon}}+\delta_{\sqrt{1-\epsilon}}\right)\]
and
\[g_\epsilon(x)=\int_\mathbb{R}\frac{s^2+1}{(s-x)^2}\;d\rho_\epsilon(s)=\frac{\epsilon x^2+
\epsilon(1-\epsilon)}{(x^2-1+\epsilon)^2}.\] Note that
$g_\epsilon(x)=\epsilon(1-\epsilon)^{-1}$ at points
$x=0,\pm\sqrt{3(1-\epsilon)}$ which implies that for
$1<t<\epsilon^{-1}$ the set $V_t^+=\{x:g_\epsilon(x)>(t-1)^{-1}\}$
consists of two components which are contained in intervals
$(-\sqrt{3(1-\epsilon)},0)$ and $(0,\sqrt{3(1-\epsilon)})$,
respectively. Moreover, when $t=\epsilon^{-1}$ the closure of the
union of these two components is the interval
$[-\sqrt{3(1-\epsilon)},\sqrt{3(1-\epsilon)}]$. Hence the support
$\mathrm{supp}((\mu_\epsilon^{\boxplus t})^{\mathrm{ac}})$ consists
of two components for $1<t<\epsilon^{-1}$ and these two components
merge into one piece when $t\geq\epsilon^{-1}$. On the other hand,
if $\epsilon<2/3$ then the measure $\mu_\epsilon^{\boxplus t}$ has
three atoms $0,\pm t$ for $1<t<2(2-\epsilon)^{-1}$, has one atom $0$
for $2(2-\epsilon)^{-1}\leq t<\epsilon^{-1}$, and has no atom for
$t\geq\epsilon^{-1}$. If $\epsilon=2/3$ then the measure
$\mu_\epsilon^{\boxplus t}$ has three atoms $0,\pm t$ for $1<t<3/2$
and has no atom for $t\geq3/2$. If $2/3<\epsilon$ then the measure
$\mu_\epsilon^{\boxplus t}$ has three atoms $0,\pm t$ for
$1<t<\epsilon^{-1}$, has two atoms $\pm t$ for $\epsilon^{-1}\leq
t<2(2-\epsilon)^{-1}$, and has no atom for
$t\geq2(2-\epsilon)^{-1}$. By Proposition 4.7 we see that
\[0\in\mathrm{supp}((\mu_\epsilon^{\boxplus
t})^{\mathrm{ac}})\;\;\;\;\;\mathrm{for}\;\;\;\;\;t\geq\epsilon^{-1}\]
and
\[\pm t\in\mathrm{supp}((\mu_\epsilon^{\boxplus
t})^{\mathrm{ac}})\;\;\;\;\;\mathrm{for}\;\;\;\;\;t\geq2(2-\epsilon)^{-1}.\]
From the above discussions, we conclude that the numbers of
components in $\mathrm{supp}(\mu_\epsilon^{\boxplus t})$ as follows.
\\ If $0<\epsilon<2/3$ then
\[n(t)=\left\{
\begin{array}{lll}
5       ,              &\hbox{for\;\;$t\in(1,2(2-\epsilon)^{-1})$} \\
3       ,              &\hbox{for\;\;$t\in[2(2-\epsilon)^{-1},\epsilon^{-1})$} \\
1 ,                    &\hbox{for\;\;$t\in[\epsilon^{-1},\infty)$}. \\
\end{array}
\right.\] If $\epsilon=2/3$ then
\[n(t)=\left\{
\begin{array}{ll}
5       ,              &\hbox{for\;\;$t\in(1,3/2)$} \\
1       ,              &\hbox{for\;\;$t\in[3/2,\infty)$}. \\
\end{array}
\right.\] Similarly, if $2/3<\epsilon<1$ then $n(t)$ begins with
$5$, reduces to $3$, and then becomes $1$ for large $t$.
\\ $\mathbf{Acknowledgments}$ The author wishes to thank his
advisor, Professor Hari Bercovici, for his generosity,
encouragement, and invaluable discussion during the course of the
investigation.

\end{document}